\input epsf
\input amssym

\input sl5.intref


\newfam\scrfam
\batchmode\font\tenscr=rsfs10 \errorstopmode
\ifx\tenscr\nullfont
        \message{rsfs script font not available. Replacing with calligraphic.}
        \def\scr{\cal}
\else   
        \font\sevenscr=rsfs7
        \font\fivescr=rsfs5
        \skewchar\tenscr='177 \skewchar\sevenscr='177 \skewchar\fivescr='177
        \textfont\scrfam=\tenscr \scriptfont\scrfam=\sevenscr
        \scriptscriptfont\scrfam=\fivescr
        \def\scr{\fam\scrfam}
        \def\cal{\scr}
\fi
\catcode`\@=11
\newfam\frakfam
\batchmode\font\tenfrak=eufm10 \errorstopmode
\ifx\tenfrak\nullfont
        \message{eufm font not available. Replacing with italic.}
        \def\frak{\it}
\else
	
	\font\sevenfrak=eufm7 \font\fivefrak=eufm5
        \font\eightfrak=eufm8
	\textfont\frakfam=\tenfrak
	\scriptfont\frakfam=\sevenfrak \scriptscriptfont\frakfam=\fivefrak
	\def\frak{\fam\frakfam}
\fi
\catcode`\@=\active
\newfam\msbfam
\batchmode\font\twelvemsb=msbm10 scaled\magstep1 \errorstopmode
\ifx\twelvemsb\nullfont\def\Bbb{\bf}
        
	\font\eightbbb=cmb10 at 8pt
	\message{Blackboard bold not available. Replacing with boldface.}
\else   \catcode`\@=11
        \font\tenmsb=msbm10 \font\sevenmsb=msbm7 \font\fivemsb=msbm5
        \textfont\msbfam=\tenmsb
        \scriptfont\msbfam=\sevenmsb \scriptscriptfont\msbfam=\fivemsb
        \def\Bbb{\relax\expandafter\Bbb@}
        \def\Bbb@#1{{\Bbb@@{#1}}}
        \def\Bbb@@#1{\fam\msbfam\relax#1}
        \catcode`\@=\active
	
	\font\eightbbb=msbm8
\fi
        \font\fivemi=cmmi5
        \font\sixmi=cmmi6
        \font\eightrm=cmr8              \def\xrm{\eightrm}
        \font\eightbf=cmbx8             \def\xbf{\eightbf}
        \font\eightit=cmti10 at 8pt     \def\xit{\eightit}

        \font\eighttt=cmtt8

        \font\eightcp=cmcsc8    
                      \def\xold{\eighti}
        \font\eightmi=cmmi8
                     \def\xbold{\eightib}
        \font\teni=cmmi10               \def\old{\teni}
        \font\tencp=cmcsc10

        \font\twelvecp=cmcsc10 scaled\magstep1
        
        \font\sixrm=cmr6
        \font\fiverm=cmr5

        \font\eightsy=cmsy8
        \font\sixsy=cmsy6
        \font\eightsl=cmsl8

        \font\sixbf=cmbx6

	 at10pt	
	\font\twelvehelvbold=phvb at12pt
	 at14pt
	\font\sixteenhelvbold=phvb at16pt
	 at16pt



\def\xbold{\xbf}
\def\xold{\xrm}


\def\noblackbox{\overfullrule=0pt}
\noblackbox

\def\eightpoint{
\def\rm{\fam0\eightrm}
\textfont0=\eightrm \scriptfont0=\sixrm \scriptscriptfont0=\fiverm
\textfont1=\eightmi  \scriptfont1=\sixmi  \scriptscriptfont1=\fivemi
\textfont2=\eightsy \scriptfont2=\sixsy \scriptscriptfont2=\fivesy
\textfont3=\tenex   \scriptfont3=\tenex \scriptscriptfont3=\tenex
\textfont\itfam=\eightit \def\it{\fam\itfam\eightit}
\textfont\slfam=\eightsl \def\sl{\fam\slfam\eightsl}
\textfont\ttfam=\eighttt \def\tt{\fam\ttfam\eighttt}
\textfont\bffam=\eightbf \scriptfont\bffam=\sixbf 
                         \scriptscriptfont\bffam=\fivebf
                         \def\bf{\fam\bffam\eightbf}
\normalbaselineskip=10pt}



\newtoks\headtext
\headline={\ifnum\pageno=1\hfill\else
	\ifodd\pageno
        \noindent{\eightcp\the\headtext}{ }\dotfill{ }{\old\folio}
	\else\noindent{\old\folio}{ }\dotfill{ }{\eightcp\the\headtext}\fi
	\fi}
\def\makeheadline{\vbox to 0pt{\vss\noindent\the\headline\break
\hbox to\hsize{\hfill}}
        \vskip2\baselineskip}
\newcount\infootnote
\infootnote=0
\newcount\footnotecount
\footnotecount=1
\def\foot#1{\infootnote=1
\footnote{${}^{\the\footnotecount}$}{\vtop{\baselineskip=.75\baselineskip
\advance\hsize by
-\parindent{\eightpoint\rm\hskip-\parindent
#1}\hfill\vskip\parskip}}\infootnote=0\global\advance\footnotecount by
1\hskip2pt}
\newcount\refcount
\refcount=1
\newwrite\refwrite
\def\oldsize{\ifnum\infootnote=1\xold\else\old\fi}
\def\ref#1#2{
	\def#1{{{\oldsize\the\refcount}}\ifnum\the\refcount=1\immediate\openout\refwrite=\jobname.refs\fi\immediate\write\refwrite{\item{[{\xold\the\refcount}]} 
	#2\hfill\par\vskip-2pt}\xdef#1{{\noexpand\oldsize\the\refcount}}\global\advance\refcount by 1}
	}
\def\refout{\eightpoint\catcode`\@=11
        \xrm\immediate\closeout\refwrite
        \vskip2\baselineskip
        \immediate\write\contentswrite{\item{}\hbox to\contentlength{References\dotfill\the\pageno}}
        {\noindent\twelvecp References}\hfill\vskip\baselineskip
        \baselineskip=.75\baselineskip
        \input\jobname.refs
        \baselineskip=4\baselineskip \divide\baselineskip by 3
        \catcode`\@=\active\rm}

\def\skipref#1{\hbox to15pt{\phantom{#1}\hfill}\hskip-15pt}

\def\hepth#1{\href{http://xxx.lanl.gov/abs/hep-th/#1}{arXiv:\allowbreak
hep-th\slash{\xold#1}}}

\def\arxiv#1#2{\href{http://arxiv.org/abs/#1.#2}{arXiv:\allowbreak
{\xold#1}.{\xold#2}}} 
 
\def\jhep#1#2#3#4{\href{http://jhep.sissa.it/stdsearch?paper=#2\%28#3\%29#4}{J. High Energy Phys. {\xbold #1#2} ({\xold#3}) {\xold#4}}}

\def\CMP#1#2#3{Commun. Math. Phys. {\xbold#1} ({\xold#2}) {\xold#3}}

\def\FP#1#2#3{Fortsch. Phys. {\xbold#1} ({\xold#2}) {\xold#3}}

\def\JMP#1#2#3{J. Math. Phys. {\xbold#1} ({\xold#2}) {\xold#3}}
\def\JPA#1#2#3{J. Phys. {\xbf A}{\xbold#1} ({\xold#2}) {\xold#3}}
\def\LMP#1#2#3{Lett. Math. Phys. {\xbold#1} ({\xold#2}) {\xold#3}}
\def\MPLA#1#2#3{Mod. Phys. Lett. {\xbf A}{\xbold#1} ({\xold#2}) {\xold#3}}

\def\NPB#1#2#3{Nucl. Phys. {\xbf B}{\xbold#1} ({\xold#2}) {\xold#3}}

\def\PLB#1#2#3{Phys. Lett. {\xbf B}{\xbold#1} ({\xold#2}) {\xold#3}}
\def\PM#1#2#3{Progr. Math. {\xbold#1} ({\xold#2}) {\xold#3}}

\newcount\sectioncount
\sectioncount=0
\def\section#1#2{\global\eqcount=0
	\global\subsectioncount=0
        \global\advance\sectioncount by 1
	\ifnum\sectioncount>1
	        \vskip2\baselineskip
	\fi
\noindent{\twelvecp\the\sectioncount. #2}\par\nobreak
       \vskip.5\baselineskip\noindent
        \xdef#1{{\old\the\sectioncount}}}
\newcount\subsectioncount
\def\subsection#1#2{\global\advance\subsectioncount by 1
\vskip.75\baselineskip\noindent\line{\tencp\the\sectioncount.\the\subsectioncount. #2\hfill}\nobreak\vskip.4\baselineskip\nobreak\noindent\xdef#1{{\old\the\sectioncount}.{\old\the\subsectioncount}}}
\def\immediatesubsection#1#2{\global\advance\subsectioncount by 1
\vskip-\baselineskip\noindent
\line{\tencp\the\sectioncount.\the\subsectioncount. #2\hfill}
	\vskip.5\baselineskip\noindent
	\xdef#1{{\old\the\sectioncount}.{\old\the\subsectioncount}}}
\newcount\subsubsectioncount
\def\subsubsection#1#2{\global\advance\subsubsectioncount by 1
\vskip.75\baselineskip\noindent\line{\tencp\the\sectioncount.\the\subsectioncount.\the\subsubsectioncount. #2\hfill}\nobreak\vskip.4\baselineskip\nobreak\noindent\xdef#1{{\old\the\sectioncount}.{\old\the\subsectioncount}.{\old\the\subsubsectioncount}}}
\newcount\appendixcount
\appendixcount=0
\def\appendix#1{\global\eqcount=0
        \global\advance\appendixcount by 1
        \vskip2\baselineskip\noindent
        \ifnum\the\appendixcount=1
        {\twelvecp Appendix A: #1}\par\nobreak
                        \vskip.5\baselineskip\noindent\fi
        \ifnum\the\appendixcount=2
        {\twelvecp Appendix B: #1}\par\nobreak
                        \vskip.5\baselineskip\noindent\fi
        \ifnum\the\appendixcount=3
        {\twelvecp Appendix C: #1}\par\nobreak
                        \vskip.5\baselineskip\noindent\fi}
\def\acknowledgements{\immediate\write\contentswrite{\item{}\hbox
        to\contentlength{Acknowledgements\dotfill\the\pageno}}
        \vskip2\baselineskip\noindent
        \underbar{\it Acknowledgements:}\ }
\newcount\eqcount
\eqcount=0
\def\Eqn#1{\global\advance\eqcount by 1
\ifnum\the\sectioncount=0
	\xdef#1{{\noexpand\oldsize\the\eqcount}}
	\eqno({\oldstyle\the\eqcount})
\else
        \ifnum\the\appendixcount=0
\xdef#1{{\noexpand\oldsize\the\sectioncount}.{\noexpand\oldsize\the\eqcount}}
                \eqno({\oldstyle\the\sectioncount}.{\oldstyle\the\eqcount})\fi
        \ifnum\the\appendixcount=1
	        \xdef#1{{\noexpand\oldstyle A}.{\noexpand\oldstyle\the\eqcount}}
                \eqno({\oldstyle A}.{\oldstyle\the\eqcount})\fi
        \ifnum\the\appendixcount=2
	        \xdef#1{{\noexpand\oldstyle B}.{\noexpand\oldstyle\the\eqcount}}
                \eqno({\oldstyle B}.{\oldstyle\the\eqcount})\fi
        \ifnum\the\appendixcount=3
	        \xdef#1{{\noexpand\oldstyle C}.{\noexpand\oldstyle\the\eqcount}}
                \eqno({\oldstyle C}.{\oldstyle\the\eqcount})\fi
\fi}
\def\eqn{\global\advance\eqcount by 1
\ifnum\the\sectioncount=0
	\eqno({\oldstyle\the\eqcount})
\else
        \ifnum\the\appendixcount=0
                \eqno({\oldstyle\the\sectioncount}.{\oldstyle\the\eqcount})\fi
        \ifnum\the\appendixcount=1
                \eqno({\oldstyle A}.{\oldstyle\the\eqcount})\fi
        \ifnum\the\appendixcount=2
                \eqno({\oldstyle B}.{\oldstyle\the\eqcount})\fi
        \ifnum\the\appendixcount=3
                \eqno({\oldstyle C}.{\oldstyle\the\eqcount})\fi
\fi}
\def\multi{\global\advance\eqcount by 1}
\def\multieqn#1{({\oldstyle\the\sectioncount}.{\oldstyle\the\eqcount}\hbox{#1})}
\def\multiEqn#1#2{\xdef#1{{\oldstyle\the\sectioncount}.{\old\the\eqcount}#2}
        ({\oldstyle\the\sectioncount}.{\oldstyle\the\eqcount}\hbox{#2})}
\def\multiEqnAll#1{\xdef#1{{\oldstyle\the\sectioncount}.{\old\the\eqcount}}}
\newcount\tablecount
\tablecount=0
\def\Table#1#2#3{\global\advance\tablecount by 1
\immediate\write\intrefwrite{\def\noexpand#1{{\noexpand\oldsize\the\tablecount}}}
       \vtop{\vskip2\parskip
       \centerline{#2}
       \vskip5\parskip
       {\narrower\noindent\it Table \the\tablecount: #3\smallskip}
       \vskip2\parskip}}
\newcount\figurecount
\figurecount=0
\def\Figure#1#2#3{\global\advance\figurecount by 1
\immediate\write\intrefwrite{\def\noexpand#1{{\noexpand\oldsize\the\figurecount}}}
       \vtop{\vskip2\parskip
       \centerline{#2}
       \vskip4\parskip
       \centerline{\it Figure \the\figurecount: #3}
       \vskip3\parskip}}
\def\TextFigure#1#2#3{\global\advance\figurecount by 1
\immediate\write\intrefwrite{\def\noexpand#1{{\noexpand\oldsize\the\figurecount}}}
       \vtop{\vskip2\parskip
       \centerline{#2}
       \vskip4\parskip
       {\narrower\noindent\it Figure \the\figurecount: #3\smallskip}
       \vskip3\parskip}}
\newtoks\url
\def\Href#1#2{\catcode`\#=12\url={#1}\catcode`\#=\active#2}
\def\href#1#2{{#2}}

\parskip=3.5pt plus .3pt minus .3pt
\baselineskip=14pt plus .1pt minus .05pt
\lineskip=.5pt plus .05pt minus .05pt
\lineskiplimit=.5pt
\abovedisplayskip=18pt plus 4pt minus 2pt
\belowdisplayskip=\abovedisplayskip
\hsize=14cm
\vsize=19cm
\hoffset=1.5cm
\voffset=1.8cm
\frenchspacing
\footline={}
\raggedbottom

\newskip\origparindent
\origparindent=\parindent

\def\*{\partial}

\def\={\!=\!}
\def\small#1{{\hbox{$#1$}}}

\def\Fraction#1#2{\small{#1\over#2}}
\def\Fr{\Fraction}

\def\eg{{\it e.g.}}

\def\ie{{\it i.e.}}



\def\appendix#1#2{\global\eqcount=0
        \global\advance\appendixcount by 1
        \vskip2\baselineskip\noindent
        \ifnum\the\appendixcount=1
        \immediate\write\intrefwrite{\def\noexpand#1{A}}
        {\twelvecp Appendix A: #2}\par\nobreak
                        \vskip.5\baselineskip\noindent\fi
        \ifnum\the\appendixcount=2
        {\twelvecp Appendix B: #2}\par\nobreak
                        \vskip.5\baselineskip\noindent\fi
        \ifnum\the\appendixcount=3
        {\twelvecp Appendix C: #2}\par\nobreak
                        \vskip.5\baselineskip\noindent\fi}

\def\textfrac#1#2{\raise .45ex\hbox{\the\scriptfont0 #1}\nobreak\hskip-1pt/\hskip-1pt\hbox{\the\scriptfont0 #2}}


\def\frac{\Fr}

\def\mathbb{\Bbb}



\def\fg{{\frak g}}

\def\sl{{\frak sl}}


\catcode`@=11
\def\openupnormal{\afterassignment\@penupnormal\dimen@=}
\def\@penupnormal{\advance\normallineskip\dimen@
  \advance\normalbaselineskip\dimen@
  \advance\normallineskiplimit\dimen@}
\catcode`@=12

\def\EqMatrix{\let\quad\enspace\openupnormal6pt\matrix}



\def\textfrac#1#2{\raise .45ex\hbox{\the\scriptfont0 #1}\nobreak\hskip-1pt/\hskip-1pt\hbox{\the\scriptfont0 #2}}


\def\frac{\Fr}

\def\mathbb{\Bbb}

\newskip\scrskip
\scrskip=-.6pt plus 0pt minus .1pt


\newwrite\intrefwrite
\immediate\openout\intrefwrite=sl5.intref

\newwrite\contentswrite

\newdimen\sublength
\sublength=\hsize 
\advance\sublength by -\parindent

\newdimen\contentlength
\contentlength=\sublength

\advance\sublength by -\parindent

\def\section#1#2{\global\eqcount=0
	\global\subsectioncount=0
        \global\advance\sectioncount by 1
\ifnum\the\sectioncount=1\immediate\openout\contentswrite=sl5.contents\fi
\immediate\write\contentswrite{\item{\the\sectioncount.}\hbox to\contentlength{#2\dotfill\the\pageno}}
        \ifnum\sectioncount>1
		        \vskip2\baselineskip
	\fi
\immediate\write\intrefwrite{\def\noexpand#1{{\noexpand\oldsize\the\sectioncount}}}\noindent{\twelvecp\the\sectioncount. #2}\par\nobreak
       \vskip.5\baselineskip\noindent}

\def\subsection#1#2{\global\advance\subsectioncount by 1
\immediate\write\contentswrite{\itemitem{\the\sectioncount.\the\subsectioncount.}\hbox
to\sublength{#2\dotfill\the\pageno}}
\immediate\write\intrefwrite{\def\noexpand#1{{\noexpand\oldsize\the\sectioncount}.{\noexpand\oldsize\the\subsectioncount}}}\vskip.75\baselineskip\noindent\line{\tencp\the\sectioncount.\the\subsectioncount. #2\hfill}\nobreak\vskip.4\baselineskip\nobreak\noindent}

\def\immediatesubsection#1#2{\global\advance\subsectioncount by 1
\immediate\write\contentswrite{\itemitem{\the\sectioncount.\the\subsectioncount.}\hbox
to\sublength{#2\dotfill\the\pageno}}
\immediate\write\intrefwrite{\def\noexpand#1{{\noexpand\oldsize\the\sectioncount}.{\noexpand\oldsize\the\subsectioncount}}}
\vskip-\baselineskip\noindent
\line{\tencp\the\sectioncount.\the\subsectioncount. #2\hfill}
	\vskip.5\baselineskip\noindent}


\def\apriori{{\it a priori}}


\def\BB{{\cal B}}


\ref\PureSGI{M. Cederwall, {\xit ``Towards a manifestly supersymmetric
    action for D=11 supergravity''}, \jhep{10}{01}{2010}{117}
    [\arxiv{0912}{1814}].}  

\ref\PureSGII{M. Cederwall, 
{\xit ``D=11 supergravity with manifest supersymmetry''},
    \MPLA{25}{2010}{3201} [\arxiv{1001}{0112}].}

\ref\CederwallKarlssonBI{M. Cederwall and A. Karlsson, {\xit ``Pure
spinor superfields and Born--Infeld theory''},
\jhep{11}{11}{2011}{134} [\arxiv{1109}{0809}].}

\ref\PureSpinorOverview{M. Cederwall, {\xit ``Pure spinor superfields
--- an overview''}, Springer Proc. Phys. {\xbf153} ({\xrm2013}) {\xrm61} 
[\arxiv{1307}{1762}].}

\ref\CederwallTensorAction{M. Cederwall, unpublished.}

\ref\CederwallBLG{M. Cederwall, {\xit ``N=8 superfield formulation of
the Bagger--Lambert--Gustavsson model''}, \jhep{08}{09}{2008}{116}
[\arxiv{0808}{3242}].}

\ref\CederwallABJM{M. Cederwall, {\xit ``Superfield actions for N=8 
and N=6 conformal theories in three dimensions''},
\jhep{08}{10}{2008}{70}
[\arxiv{0809}{0318}].}

\ref\CederwallReformulation{M. Cederwall, {\xit ``An off-shell superspace
reformulation of $D=4$, $N=4$ super-Yang--Mills theory''},
\FP{66}{2018}{1700082} [\arxiv{1707}{00554}].}

\ref\CederwallDSix{M. Cederwall, {\xit ``Pure spinor superspace action
for D=6, N=1 super-Yang--Mills theory''}, \jhep{18}{05}{2018}{115}
[\arxiv{1712.02284}]} 

\ref\CederwallExotic{M. Cederwall {\xit ``Superspace formulation of
exotic supergravities in six dimensions''}, \jhep{21}{03}{2021}{56} [\arxiv{2012}{02719}].}

\ref\CederwallPalmkvistSaberiInProgress{M. Cederwall, J, Palmkvist and
I. Saberi, work in progress.}

\ref\CederwallPalmkvistBorcherds{M. Cederwall and J. Palmkvist, {\xit
``Superalgebras, constraints and partition functions''},
\jhep{08}{15}{2015}{36} [\arxiv{1503}{06215}].}

\ref\EagerSaberiWalcher{R. Eager, I. Saberi and J. Walcher, {\xit
``Nilpotence varieties''}, Ann. Henri Poincar\'e {\xbf22} (2021) 1319 [\arxiv{1807}{03766}].}

\ref\SpinorialCohomology{M. Cederwall, B.E.W. Nilsson and D. Tsimpis, 
{\xit ``Spinorial cohomology and maximally supersymmetric theories''},
\jhep{02}{02}{2002}{009} [\hepth{0110069}].}

\ref\BerkovitsI{N. Berkovits, 
{\xit ``Super-Poincar\'e covariant quantization of the superstring''}, 
\jhep{00}{04}{2000}{018} [\hepth{0001035}].}

\ref\BerkovitsIII{N. Berkovits, 
{\xit ``Cohomology in the pure spinor formalism for the
superstring''}, 
\jhep{00}{09}{2000}{046} [\hepth{0006003}].}

\ref\Movshev{M. Movshev and A. Schwarz, {\xit ``On maximally
supersymmetric Yang--Mills theories''}, \NPB{681}{2004}{324}
[\hepth{0311132}].}

\ref\MovshevSchwarzDef{M. Movshev and A. Schwarz, {\xit
``Supersymmetric deformations of maximally supersymmetric gauge
theories''}, \jhep{12}{09}{2012}{136} [\arxiv{0910}{0620}].}

\ref\CederwallPalmkvistExtendedGeometry{M. Cederwall and J. Palmkvist,
{\xit ``Extended geometries''}, \jhep{02}{18}{2018}{071} [\arxiv{1711}{07694}].}

\ref\CederwallPalmkvistLInfty{M. Cederwall and J. Palmkvist, {\xit
``$L_\infty$ algebras for extended geometry from Borcherds
superalgebras''}, \CMP{369}{2019}{2} [\arxiv{1804}{04337}].}

\ref\MovshevSchwarzAlgebra{M. Movshev and A. Schwarz, {\xit
``Algebraic structure of Yang--Mills theory''}, \PM{244}{2006}{473}
[\hepth{0404183}].} 

\ref\MovshevDeform{M. Movshev, {\xit ``On deformations of Yang--Mills
algebras''}, \hepth{0509119}.}

\ref\PalmkvistTensor{J. Palmkvist, {\xit ``The tensor hierarchy
algebra''}, \JMP{55}{2014}{011701} [\arxiv{1305}{0018}].}

\ref\PalmkvistHierarchy{J. Palmkvist, {\xit ``Tensor hierarchies,
Borcherds algebras and $E_{11}$''}, \jhep{12}{02}{2012}{066}
[\arxiv{1110}{4892}].} 

\ref\PalmkvistExceptionalGeometry{J. Palmkvist, {\xit ``Exceptional
geometry and Borcherds superalgebras''}, \jhep{15}{11}{2015}{032} [\arxiv{1507}{08828}].}

\ref\CarboneCederwallPalmkvist{L. Carbone, M. Cederwall and
J. Palmkvist, {\xit ``Generators and relations for Lie superalgebras
of Cartan type''}, \JPA{52}{2019}{055203} [\arxiv{1802}{05767}].}

\ref\CederwallPalmkvistTHAI{M. Cederwall and J. Palmkvist, {\xit
``Tensor hierarchy algebras and extended geometry I: Construction of
the algebra''}, \jhep{20}{02}{2020}{144} [\arxiv{1908}{08695}].}

\ref\CederwallPalmkvistTHAII{M. Cederwall and J. Palmkvist, {\xit
``Tensor hierarchy algebras and extended geometry II: Gauge structure
and dynamics''}, \jhep{20}{02}{2020}{145} [\arxiv{1908}{08696}].}

\ref\CederwallPalmkvistHyperbolicTHA{M. Cederwall and J. Palmkvist,
{\xit ``Tensor hierarchy extensions of hyperbolic Kac--Moody
algebras''}, \arxiv{2103}{02476}.}

\ref\BeyondEEleven{G. Bossard, A. Kleinschmidt, J. Palmkvist,
C.N. Pope and E. Sezgin, {\xit ``Beyond $E_{11}$''},
\jhep{17}{05}{2017}{020} [\arxiv{1703}{01305}].}

\ref\KacRudakov{V.G. Kac and A. Rudakov, {\xit ``Complexes of modules
over exceptional Lie superalgebras $E(3,8)$ and $E(5,10)$''},
Int. Math. Res. Notices {\xbf19} (2002) 1007
[arXiv:math-ph/{\xold0112022}].}

\ref\ChengKac{S.-J. Cheng and V.G. Kac, {\xit ``Structure of some
{\eightbbb Z}-graded Lie superalgebras of vector fields''},
Transformation Groups {\xbf4} (1999) 219.}

\ref\KacLinComp{V.G. Kac, {\xit ``Classification of
infinite-dimensional simple linearly compact Lie superalgebras''},
Adv. Math. {\xbf139} (1998) 1.}

\ref\CantariniKac{N. Cantarini and V.G. Kac, {\xit ``Automorphisms and
forms of simple infinite-dimensional linearly compact Lie
superalgebras''}, Int. J. Geom. Meth. Mod. Phys. {\xbf3} (2006) 845
[arXiv:math/0601292].}

\ref\CantariniCaselli{N. Cantarini and F. Caselli, {``\xit Low degree
morphisms of $E(5,10)$-generalized Verma modules''}, Alg. Repr. Theory
{\xbf23} (2020) 2131 [\arxiv{1903}{11438}].}

\ref\Shchepochkina{I.M. Shchepochkina, {\xit ``Five exceptional simple
	Lie superalgebras of vector fields''}, Functional Analysis and
	its Applications {\xbf33} (1999) 208.}

\ref\BerkovitsGuillen{N. Berkovits and M. Guillen, {\xit ``Equations
    of motion from Cederwall's pure spinor superspace
    actions''}, \jhep{18}{08}{2018}{033} [\arxiv{1804}{06979}].}

\ref\JonssonMSc{S. Jonsson, MSc thesis, Chalmers Univ. of Technology ({\xold2021}).}

\ref\PalmkvistPrivate{J. Palmkvist, private communication.}

\ref\KacSuperalgebrasII{V.G. Kac, {\xit ``Lie superalgebras''},
Adv. Math. {\xbf26} (\xold{1977}) {\xold8}.}

\ref\ChestermanGhost{M. Chesterman, {\xit ``Ghost constraints and the
covariant quantisation of the superparticle in ten dimensions''},
\jhep{04}{02}{2004}{011}, [\hepth{0212261}].}

\ref\BerkovitsParticle{N. Berkovits, {\xit ``Covariant quantization of
the superparticle using pure spinors''}, \jhep{01}{09}{2001}{016}
[\hepth{0105050}].}

\ref\BerkovitsNekrasovCharacter{N. Berkovits and N. Nekrasov, {\xit
    ``The character of pure spinors''}, \LMP{74}{2005}{75}
  [\hepth{0503075}].}

\ref\BermanCederwallKleinschmidtThompson{D.S. Berman, M. Cederwall,
A. Kleinschmidt and D.C. Thompson, {\xit ``The gauge structure of
generalised diffeomorphisms''}, \jhep{13}{01}{2013}{64} [\arxiv{1208}{5884}].}

\ref\SaberiWilliams{I. Saberi and B. Williams, {\xit ``Twisting pure
spinor superfields, with applications to supergravity''},
\arxiv{2106}{15639}.}

\ref\BaulieuSUFive{L. Baulieu, {\xit ``SU(5)-invariant decomposition
of ten-dimensional Yang-Mills supersymmetry''},
\PLB{698}{2011}{63} [\arxiv{1009}{3893}].}


%
\line{
\epsfysize=18mm
\epsffile{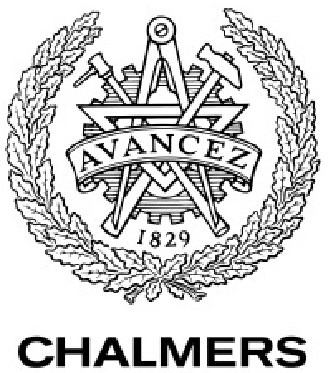}
\hfill}
\vskip-16mm

\line{\hfill}
\line{\hfill Gothenburg preprint}
\line{\hfill July, {\old2021}}
\line{\hrulefill}


\headtext={Cederwall: 
``SL(5) supersymmetry''}

\vfill

\centerline{\sixteenhelvbold
SL(5) supersymmetry}


%

\vfill

\centerline{\twelvehelvbold Martin Cederwall}

\vfill
\vskip-1cm

\centerline{\it Department of Physics}
\centerline{\it Chalmers University of Technology}
\centerline{\it SE 412 96 Gothenburg, Sweden}

\vfill

{\narrower\noindent \underbar{Abstract:}
We consider supersymmetry in five dimensions, where the fermionic
parameters are a $2$-form under $SL(5)$. Supermultiplets are
investigated using the pure spinor superfield formalism, and are found
to be closely related to infinite-dimensional extensions of the
supersymmetry algebra: the Borcherds superalgebra $\BB(E_4)$,
the tensor hierarchy algebra $S(E_4)$ and the
exceptional superalgebra $E(5,10)$. A theorem relating $\BB(E_4)$ and
$E(5,10)$ to all levels is given.
\smallskip}
\vfill

\font\xxtt=cmtt6

\vtop{\baselineskip=.6\baselineskip\xxtt
\line{\hrulefill}
\catcode`\@=11
\line{email: martin.cederwall@chalmers.se\hfill}
\catcode`\@=\active
}

\eject


\catcode`\@=11
        \vskip2\baselineskip
        {\noindent\twelvecp Contents}\hfill\vskip\baselineskip
        \input sl5.contents
        \catcode`\@=\active\rm
\vskip3\baselineskip

\section\IntroSection{Introduction and overview}Supersymmetry provides an
extension of bosonic space-time symmetries with fermionic
generators. These are generically spinors under space-time rotations
(and may also transform under R-symmetry). In certain situations,
supersymmetry generators in non-spinorial modules may be considered.
The main example is provided by ``twisting'', where one considers a
fermionic generator which is a singlet under some subalgebra.

More generically, one may \apriori\ consider an assignment where
supersymmetry generators come in a module $S$ of a space-time
``structure group'' $G$, which we think of as corresponding to the double
cover of the Lorentz group together with R-symmetry. 
A supersymmetry algebra will take the form\foot{We use a
notation with $[\cdot,\cdot]$ for the graded Lie brackets or graded
commutators; the present bracket is of course symmetric.}
$[Q_a,Q_b]=c_{ab}{}^mP_m$, with some invariant tensor $c$, and the
rest of the brackets vanishing.
The only condition is that the symmetric
product $\vee^2S$ contains the vector representation $V$.

Presently, we will consider one specific such assignment, namely when
the structure group is $G=SL(5)$, $V=\overline{\bf5}$ and $S={\bf10}$.
The supersymmetry algebra then is\foot{A factor $i$ may be included in
the right hand side, depending on conventions.}
$$
[Q^{mn},Q^{pq}]=2\epsilon^{mnpqr}\partial_r\;,\Eqn\SLFiveSSAlgebra
$$
There are several observations that make this choice special, and
exploring them is the purpose of the paper.

A general method for formulating supersymmetric field theories on
superspace is provided by ``pure spinor superfield theory''
[\SpinorialCohomology\skipref\BerkovitsI\skipref\BerkovitsParticle\skipref\BerkovitsIII\skipref\PureSGI\skipref\PureSGII\skipref\CederwallKarlssonBI\skipref\CederwallBLG\skipref\CederwallABJM\skipref\CederwallReformulation\skipref\CederwallDSix\skipref\CederwallExotic\skipref\BerkovitsGuillen-\PureSpinorOverview], where the
superspace coordinates $x\in V$ and $\theta\in S$ are complemented by
a ``bosonic spinor'' $\lambda\in S$, which is subject to the
constraint $c_{ab}{}^m\lambda^a\lambda^b=0$.
The content of supermultiplets can be deduced already from the
partition function (Hilbert series) of $\lambda$.

In the simple case of $D=10$ super-Yang--Mills theory, the
constraint on $\lambda$ turns it into a Cartan pure spinor. This
(more generally, the fact that $\lambda$ belongs to a minimal
$S$-orbit under $G$) 
enables a Koszul duality [\CederwallPalmkvistBorcherds] between the
associative algebra generated by $\lambda$ and the positive levels of
a Borcherds superalgebra. The structure of the Borcherds superalgebra,
in the $D=10$ super-Yang--Mills case $\BB(E_5)$,
is thus closely connected to the supersymmetry multiplet
[\MovshevSchwarzAlgebra,\MovshevDeform,\CederwallPalmkvistSaberiInProgress].

In  Section \KoszulSection, we will establish such a relation for the $SL(5)$
supersymmetry algebra (\SLFiveSSAlgebra), the Borcherds superalgebra
$\BB(E_4)$ and a certain supermultiplet\foot{The
notation $E_4\simeq A_4$ is used since it is part of the $E$ series,
and the super-extension is associated to the leftmost node in the
Dynkin diagram, see Figure \BEFourFigure.}.
This is achieved by applying, in Section \ScalarFieldCohomologySection,
the principles of pure spinor superfield
theory to the case at hand.
This supermultiplet found turns out, as a vector space, to be the adjoint
module of the (infinite-dimensional) exceptional superalgebra
$E(5,10)$ [\KacLinComp,\ChengKac,\Shchepochkina], of which the global
supersymmetry algebra (\SLFiveSSAlgebra) is a subalgebra.
This observation then leads to a surprising relation
between $E(5,10)$ and the Borcherds superalgebra $\BB(E_4)$ or the
tensor hierarchy algebra [\PalmkvistTensor] $S(E_4)$.
More precisely, the half of $S(E_4)$ at levels $\geq3$ turns out to be
freely generated by the coadjoint module of $E(5,10)$.
We also describe $E(5,10)$ as a ``restricted tensor hierarchy
algebra'' in terms of the generalisation of Chevalley generators
introduced in refs. [\CarboneCederwallPalmkvist,\CederwallPalmkvistTHAI].

It may be noted that the supersymmetry algebra (\SLFiveSSAlgebra) is a
subalgebra of the $D=10$, $N=1$ supersymmetry algebra.
Under the subgroup $SL(5)\subset Spin(10)$, the vector and spinor
branch as
${\bf10}\rightarrow{\bf 5}\oplus\overline{\bf5}$,
${\bf16}\rightarrow{\bf1}\oplus{\bf 10}\oplus\overline{\bf5}$, and
the module ${\bf10}$ parametrises the infinitesimal (projective)
deformations 
of a pure spinor ${\bf1}$.
The (on-shell) $D=10$ super-Yang--Mills supermultiplet, effectively encoded in
a pure spinor of $Spin(10)$, must be possible to describe in terms of
$SL(5)$ supermultiplets. We shall comment on this in
Section \OtherSuperfieldsSection. 

Infinite-dimensional superalgebras, in particular Borcherds
superalgebras or tensor hierarchy algebras, play an important r\^ole
in the context of extended geometry
[\PalmkvistHierarchy,\PalmkvistExceptionalGeometry,\CederwallPalmkvistExtendedGeometry,\CederwallPalmkvistTHAII].
The superalgebra underlying the extended geometry for $D=11$
supergravity reduced to $d$ dimensions
is $S(E_{12-d})$. In the series, we find $S(E_4)$, which (most likely)
coincides with
$\BB(E_4)$ at positive levels, governing the gauge structure.
$S(E_5)$ (or $\BB(E_5)$) lies behind extended geometry for maximal
supergravity in 
$d=7$, but is also central to $D=10$ super-Yang--Mills theory. 
It would be interesting to understand whether there is some deeper 
reason that the exact same algebraic structures appear in two
seemingly different contexts.

The Koszul duality described in Section \KoszulSection\ holds for any
object in a minimal orbit and the corresponding Borcherds
superalgebra.
Our impression is that this is not enough for the meaningful description of a
supersymmetry multiplet, which also seems to rely on the superalgebra
having some freely generated part.
This happens also in cases, such as $D=11$ supergravity, where the
orbit is not minimal, and the dual superalgebra is not even a Lie
superalgebra [\JonssonMSc,\CederwallPalmkvistSaberiInProgress].
One may for example ask if $\BB(E_6)$ or $S(E_6)$ (which differ
slightly even at positive levels [\PalmkvistPrivate]) is relevant for
some kind of supersymmetry based on $E_6$. 

We expect relations similar to the one between $E(5,10)$ and $S(E_4)$
to hold also for some other of the exceptional superalgebras
[\KacLinComp,\ChengKac,\Shchepochkina], such as
$E(3,6)$ or $E(3,8)$, that may have a relation to $S(E_3)$, but this
issue has yet to be investigated.

\section\SLFiveSSSection{$SL(5)$ supersymmetry}The supersymmetry
algebra (\SLFiveSSAlgebra) is realised by
$$
Q^{mn}={\partial\over\partial\theta_{mn}}
+\epsilon^{mnpqr}\theta_{pq}\partial_r\;,\eqn
$$
with $\partial_m={\partial\over\partial x^m}$. Then,
$[D^{mn},Q^{pq}]=0$, where the covariant fermionic derivative is
$$
D^{mn}={\partial\over\partial\theta_{mn}}
-\epsilon^{mnpqr}\theta_{pq}\partial_r\;.\eqn
$$
Now,
$$
[D^{mn},D^{pq}]=-2\epsilon^{mnpqr}\partial_r\;,\eqn
$$
and the superspace torsion is
$T^{mn,pq,r}=2\epsilon^{mnpqr}$.

The BRST operator of pure spinor field theory is, as usual,
constructed as
$$
Q=\lambda_{mn}D^{mn}\;.\Eqn\PSFTQ
$$
Its nilpotency is guaranteed by the bilinear constraint
$$
\lambda_{[mn}\lambda_{pq]}=0\;.\eqn
$$
Denoting representations and representation modules by the Dynkin
label of the highest weight, and letting $\lambda\in S=(0100)=R(\Lambda_2)$, the
general symmetric product is
$$
\vee^2S=\vee^2(0100)=(0200)\oplus(0001)\;.
$$
The constraint implies that only the module with highest weight
$2\Lambda_2$ survives, so $\lambda$ belongs to the (unique) minimal
$S$-orbit under $SL(5)$, which is a c\^one over the $6$-dimensional Grassmannian
$Gr(2,5)$.   

\subsection\ScalarFieldCohomologySection{The cohomology of a scalar
superfield}We will be quite brief about the technicalities of the
calculations of ``pure spinor field theory'' leading to the
supermultiplets in this and the following subsection. They go along
the general principles explained \eg\ in ref. [\PureSpinorOverview].
The calculation of the zero-mode cohomologies is a matter of comparing
components in different superfields, a pure algebraic problem well
suited for a computer.
Here, we choose to put stronger focus on the way in which a
supermultiplet appears in the partition function of the constrained
object $\lambda$, and in the following Section on the Koszul duality to a superalgebra.

``Physical states'' may be defined as cohomology of the
BRST operator 
$Q$ of eq. (\PSFTQ).
They are also directly encoded in the partition function of
$\lambda$.
This partition function encodes the modules $S_p$ of monomials of
degree of homogeneity $n$ in $\lambda$ as the
coefficient of $t^n$ in a formal power series with coefficients in the
representation ring as $Z_\lambda(t)=\oplus_{p=0}^\infty S_pt^p$. We
choose the conventions that $S_p={R(p\Lambda_3)}$
are the modules of the {\it
components} in the expansion\foot{The duality
formulated in ref. [\CederwallPalmkvistBorcherds] and in
Section \KoszulSection\ employs a relation to the {\it coalgebra}. We
account for this by this definition of the partition function; an
alternative would be to employ negative instead of positive levels.}, which are conjugate to
the ones of the basis elements
$\lambda_{m_1n_1}\ldots\lambda_{m_pn_p}$.
Thus, 
$$
Z_\lambda(t)=\bigoplus_{p=0}^\infty(00p0)t^p\;.\eqn
$$
Factoring out the partition function of an unconstrained 
object in $S$, which we denote
$(1-t)^{-(0010)}\equiv\oplus_{p=0}^\infty\vee^p(0010)t^p$ (and which
is compensated by $\theta$),
$$
Z_\lambda(t)=(1-t)^{-(0010)}\otimes\left((0000)\ominus(1000)t^2\oplus(0001)t^3
\ominus(0000)t^5\right)\;.\eqn
$$
The interpretation of the numerator is as the zero-mode
cohomology, \ie, the cohomology of
$\lambda_{mn}{\partial\over\partial\theta_{mn}}$ on a scalar field
$\Psi(\theta,\lambda)$. Assigning ghost number $1$ to $\Psi$ (and of
course $0$ to $\theta$ and $1$ to $\lambda$), the interpretation of
the zero-mode cohomology is:
\item{$\bullet$}A ghost $c$;
\item{$\bullet$}A $1$-form $\alpha$, appearing in $\Psi$ as
$\epsilon^{mnpqr}\lambda_{mn}\theta_{pq}\alpha_r$;
\item{$\bullet$}A vector $\xi$, appearing as
$\epsilon^{mnpqr}\lambda_{mn}\theta_{pq}\theta_{rs}\xi^s$;
\item{$\bullet$}An ``antifield'' $\gamma$, appearing as
$\epsilon^{mnrst}\epsilon^{pquvw}\lambda_{mn}\lambda_{pq}
\theta_{rs}\theta_{tu}\theta_{vw}\gamma$.

\Table\CohomologyTableScalar
{$$\hskip-2cm
\vtop{\baselineskip25pt\lineskip0pt
\ialign{
$\hfill#\quad$&$\,\hfill#\hfill\,$&$\,\hfill#\hfill\,$
&$\,\hfill#\hfill\,$&$\,\hfill#\hfill\,$\cr
&\lambda^0&\lambda^1&\lambda^2&\lambda^3\cr
&(0000)&\phantom{(0000)}&\phantom{(0000)}&\phantom{(0000)}\cr
&\bullet&\bullet            \cr 
&\bullet&(1000)&\bullet&      \cr
&\bullet&(0001)&\bullet&\bullet\cr
&\bullet&\bullet&\bullet&\bullet\cr
&\bullet&\bullet&(0000)&\bullet\cr
&\bullet&\bullet&\bullet&\bullet\cr
}}
$$}
{The zero-mode cohomology in $\Psi$. The
superfields at different ghost numbers are shifted so that fields in
the same row has the same dimension. Black dots denote the absence of
cohomology.}

Further extracting a factor
$(1-t^2)^{(1000)}\equiv\oplus_{i=0}^5\wedge^i(1000)t^{2i}$ (which
is compensated by the $x$-dependence),
$$
Z_\lambda(t)=(1-t)^{-(0010)}\otimes(1-t^2)^{(1000)}\otimes\Bigl((0000)
\oplus\bigoplus_{i=0}^\infty(i001)t^{3+2i}
\ominus\bigoplus_{i=0}^\infty(i100)t^{4+2i}
\Bigr)\;.\Eqn\FieldPartition
$$
The terms in the rightmost factor are interpreted as
the ghost zero-mode together with the derivative expansions of
$\xi$ and $\chi=d\alpha$.

In the full cohomology, closedness implies the ``equation of motion''
$\partial_m\xi^m=0$, and quotienting out exact functions gives the
gauge invariance $\alpha\sim\alpha+d\beta$. This cohomology is
reflected by the partition function (\FieldPartition).
This off-shell supermultiplet has $4$ bosonic and $4$ fermionic
local degrees of freedom. As we will see later, it corresponds to the
exceptional Lie superalgebra $E(5,10)$.
The same cohomology arises from twisting of $D=11$ supergravity
[\SaberiWilliams].

\subsection\OtherSuperfieldsSection{Cohomology of other
superfields}This Section lies outside the main line of the paper, and
may be skipped in a linear reading.

One may attempt to assign some non-trivial module to a
field. Then, some so-called shift symmetry [\CederwallKarlssonBI]
must also be introduced, otherwise the cohomology would just become the
tensor product of the module of the field with the cohomology found 
in a scalar field. In the language of ref. [\EagerSaberiWalcher] this amounts to
considering not functions on the minimal orbit, but sections of some
geometric sheaf (such as products of the (co-)tangent sheaf) over
the minimal orbit.

Taking a field $\Phi^m(x,\theta,\lambda)$ in the vector module, and
imposing equivalence under the shift symmetry
$\Phi^m\sim\Phi^m+\lambda_{np}\varrho^{mnp}$ for completely
antisymmetric $\varrho^{mnp}(x,\theta,\lambda)$ leads to the zero-mode
cohomology of Table \CohomologyTableVector.

\Table\CohomologyTableVector
{$$\hskip-2cm
\vtop{\baselineskip25pt\lineskip0pt
\ialign{
$\hfill#\quad$&$\,\hfill#\hfill\,$&$\,\hfill#\hfill\,$
&$\,\hfill#\hfill\,$&$\,\hfill#\hfill\,$\cr
&\lambda^0&\lambda^1&\lambda^2\cr
&(0001)&\phantom{(0000)}&\phantom{(0000)}\cr
&(0100)&\bullet            \cr 
&\bullet&\bullet&\bullet    \cr
&\bullet&(0010)&\bullet\cr
&\bullet&(1000)&\bullet\cr
&\bullet&\bullet&\bullet\cr
}}
$$}
{The zero-mode cohomology in a vector field $\Phi^m$.}

In the full cohomology, the vector $v^m(x)$ will obey $\*_m\*_nv^n=0$
and the $2$-form $\omega$ is closed. In this sense, the same local degrees of
freedom are reproduced as the ones in a scalar field (there is only
a single extra linear singlet mode in the derivative expansion of $v$).
On the other hand, if $\Psi$ is taken to be fermionic and of ghost
number $1$, and $\Phi^m$ bosonic of ghost number $0$, they can be
combined into a description of the $D=10$ super-Yang--Mills multiplet
with $SL(5)$ (or $SU(5)$) covariance [\BaulieuSUFive].
Then the $1$-form $\alpha$ and
the vector $v$ are the components of the $D=10$ connection, and
the vector $\xi$ and the $2$-form $\omega$ are part of the spinor.

Yet another possibility is a $1$-form field $\Xi_m(x,\theta,\lambda)$,
with shift symmetry
$\Xi_m\sim\Xi_m+\lambda_{mn}\varrho^n$. The zero-mode cohomology is
given in Table \CohomologyTableOneForm.

\Table\CohomologyTableOneForm
{$$\hskip-2cm
\vtop{\baselineskip25pt\lineskip0pt
\ialign{
$\hfill#\quad$&$\,\hfill#\hfill\,$&$\,\hfill#\hfill\,$
&$\,\hfill#\hfill\,$&$\,\hfill#\hfill\,$\cr
&\lambda^0&\lambda^1&\lambda^2\cr
&(1000)&\phantom{(0000)}&\phantom{(0000)}\cr
&(0001)&\bullet            \cr 
&\bullet&(2000)&\bullet    \cr
&\bullet&(0000)\oplus(1001)&\bullet\cr
&\bullet&(0010)&\bullet\cr
&\bullet&\bullet&\bullet\cr
}}
$$}
{The zero-mode cohomology in a $1$-form field $\Xi_m$.}

\noindent The supermultiplet at $\lambda^1$ is interesting.
There is a symmetric
tensor $h_{mn}$ with a gauge transformation
$\delta_uh_{mn}=\*_{(m}u_{n)}$, which looks like the linearised
transformation of a gravity field. The fermions $\psi_m{}^n$ have a
gauge symmetry $\delta_\epsilon\psi_m{}^n=\*_m\epsilon^n$, and there
is a bosonic $3$-form (field strength).
There are $20$ bosonic and $20$ fermionic
local degrees of freedom. The multiplet will certainly be a part of a
decomposition of $D=10$, $N=1$ supergravity, but may also have 
some significance of its own.

There may be more interesting supermultiplets that we have not found.

\section\KoszulSection{Minimal orbit partition function and Koszul
duality}\Figure\BEFourFigure{\epsffile{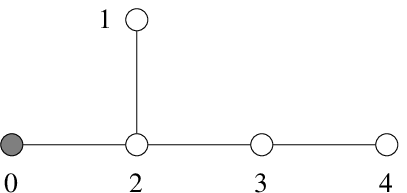}}{The Dynkin diagram of
$\BB(E_4)$, $S(E_4)$ and $E(5,10)$.}

\noindent It is known that there in general is a Koszul duality between
the associative algebra generated by an object $\lambda\in S$ in a minimal
orbit under $G$ and the positive part $\BB_+$ of a Borcherds superalgebra 
[\CederwallPalmkvistBorcherds].
The Lie superalgebra $\BB$ is defined by a Dynkin diagram obtained by
extending the Dynkin diagram of $\fg=\hbox{Lie}(G)$ by a ``fermionic''
null root connected according to the Dynkin label of $S$. The diagram
is depicted in Figure \BEFourFigure.

The concrete relation is
$$
Z_\lambda(t)\otimes Z_{\BB_+}(t)=1\;,\Eqn\KoszulRelation
$$
where $Z_{\BB_+}$ is the partition function of (the universal enveloping
algebra of) $\BB_+$, defined as
$$
Z_{\BB_+}(t)=\bigotimes_{p=1}^\infty(1-t^p)^{-(-1)^pR_p}\;,\eqn
$$
$R_p$ being the module of the level $p$ generators in $\BB$.
The duality (\KoszulRelation) can be understood as a factorised
version of the summation form
$Z_\lambda(t)=\bigoplus_{p=0}^\infty\overline{R(p\Lambda)}\,t^p$, with $\Lambda$
the highest weight of $S$, and thus as a denominator formula for
$\BB$.
A concrete way of understanding the duality is in terms of the BRST
operator for the bilinear constraint on $\lambda$. The modules $R_p$
are the modules of the corresponding infinite tower of ghosts
[\CederwallPalmkvistBorcherds,\ChestermanGhost,\BerkovitsNekrasovCharacter,\BermanCederwallKleinschmidtThompson]. This BRST operator can be identified
as the coalgebra differential of the co-superalgebra
$\BB_+^\star$. The identification relies on the absence of Lie
superalgebra cohomology other than polynomials of level $1$, which
holds for Borcherds superalgebras.

In the present case, we thus have
$$
\eqalign{
Z_{\BB_+}(t)&=(1-t)^{(0010)}\otimes(1-t^2)^{-(1000)}\cr
&\qquad\otimes\Bigl(1
\oplus\bigoplus_{i=0}^\infty(i001)t^{3+2i}
\ominus\bigoplus_{i=0}^\infty(i100)t^{4+2i}
\Bigr)^{-1}\cr
&=(1-t)^{(0010)}\otimes(1-t^2)^{-(1000)}\otimes(1-t^3)^{(0001)}
\otimes(1-t^4)^{-(0100)}\cr
&\qquad\otimes(1-t^5)^{(1001)}
\otimes(1-t^6)^{-(0002)-(1100)}\otimes\ldots
\;
}\Eqn\BEFourPartition
$$

The factor $\left(1
\oplus\bigoplus_{i=0}^\infty(i001)t^{3+2i}
\ominus\bigoplus_{i=0}^\infty(i100)t^{4+2i}
\right)^{-1}$ is now read as the partition function for the
superalgebra freely generated\foot{The partition function for an
algebra freely generated by generators in $P(t)$ is $(1-P(t))^{-1}$.
Even if the freely generated algebra
itself is complicated to describe level by level, its universal enveloping
algebra has the simple partition function
$\bigoplus_{i=0}^\infty\otimes^iP(t)=(1-P(t))^{-1}$.}
by the supermultiplet with partition
$$
P(t)=\ominus\bigoplus_{i=0}^\infty(i001)t^{3+2i}
\oplus\bigoplus_{i=0}^\infty(i100)t^{4+2i}\;.\eqn
$$
Note that this partition function is consistently graded in the sense
that when the products $\otimes^iP(t)$ are evaluated, all terms at
odd level are negative and at even level positive, so no cancellations
may arise. One can then safely identify the algebra at levels $\geq3$
as freely generated by the multiplet.
The picture is in complete analogy with how the Borcherds superalgebra
$\BB(E_5)$ is freely generated by the $D=10$ Yang--Mills on-shell
supermultiplet from level $3$
[\MovshevSchwarzAlgebra,\MovshevDeform,\CederwallPalmkvistSaberiInProgress]. 

We believe that the Koszul duality in general can be extended to fields
in non-trivial modules with
shift symmetry. Then, just as the Koszul duality described above can
be seen as providing a denominator formula for the Borcherds superalgebra, it
may be conjectured that a corresponding relation involving a
non-scalar field will provide a character formula for a representation
of the superalgebra [\CederwallPalmkvistSaberiInProgress].

\section\SuperalgebrasSection{Superalgebras based
on $SL(5)$ supersymmetry}There is a number of infinite-dimensional
superalgebras that all contain generators in $(0010)$ at level $1$ and
$(1000)$ at level $2$. They can all be described by the same Dynkin
diagram, Figure \BEFourFigure.
They have standard Chevalley generators $e_a$, $a=0,\ldots,4$ and
$f_i$, $h_i$, $i=1,\ldots,4$, but differ in terms of the level $-1$
generators. The Borcherds superalgebra $\BB(E_4)$ is
contragredient, having a Chevalley generators $f_0$ and $h_0$
with the usual relations.
The generator $e_0$, with the Serre relations $[e_0,e_0]=0$ following
from the Cartan matrix, is the lowest weight state in the $2$-form
module at level $1$, and the Serre relation implies that level $2$ only
contains a vector.


\subsection\SEFourSection{$S(E_4)$}While the superalgebra
$\BB(E_4)$ is contragredient, so the module at level $-p$ is conjugate
to the one at level $p$, this does not apply to the tensor hierarchy
algebras [\PalmkvistTensor] $S(E_4)$.
In $S(E_4)$, $h_0$ is removed, and
the generator $f_0$ of the Borcherds superalgebra is
replaced by $f_{0j}$, $i=1,3,4$ (the nodes not connected to node
$0$). New brackets are $[h_i,f_{0j}]=-A_{0i}f_{0j}$,
$[e_0,f_{0j}]=h_j$
[\CarboneCederwallPalmkvist,\CederwallPalmkvistTHAI,\CederwallPalmkvistHyperbolicTHA].
The presence of $f_{01}$ implies the presence of $(2000)$ at level
$-1$, while $f_{03},f_{04}$ give $(0011)$. The reducibility of level
$-1$ comes from the disconnectedness of the Dynkin diagram of $A_4$
with node $2$ removed.


Although it has not been proven, we strongly believe that 
the positive levels of $\BB(E_4)$ and $S(E_4)$ coincide,
$S_+(E_4)\simeq\BB_+(E_4)$,
and that
consequently the Koszul
duality involving the supersymmetry multiplet, as described in
Section \KoszulSection, applies equally well to
$S(E_4)$.
It seems likely that a proof of this may be based on an argument that the freely
generated property of the part of $\BB(E_4)$ at levels $\geq3$ does
not allow for any ideal annihilated by the negative level generators
in $S(E_4)$. This will be postponed to future examination
[\CederwallPalmkvistSaberiInProgress].

Instead of standard contragredience, $S(E_4)$ allows for an invariant
bilinear form pairing level $p$ with level $5-p$, so that
$R_{5-p}=\overline{R_p}$ [\BeyondEEleven].

\subsection\EFiveTenSection{$E(5,10)$}The superalgebra $E(5,10)$
[\KacLinComp,\ChengKac,\Shchepochkina]
is one of the ``exceptional'' simple superalgebras which are
linearly compact (which for our purposes holds if the elements arise
as a power series of the coordinates of some finite-dimensional space)
and of finite depth (meaning that there is a
maximal level\foot{In most of the mathematical literature,
level is defined with a minus sign compared to our conventions, so 
this would
read as minimal level.}), and has attracted some interest in the mathematics
literature [\KacRudakov,\CantariniKac,\CantariniCaselli].
It is defined as a super-extension of volume-preserving diffeomorphisms in
$5$ dimensions by fermionic generators, the parameters of which are
closed $2$-forms. Letting $\xi$, $\eta$ be divergence-free vector
fields and $\chi$, $\psi$ closed fermionic $2$-forms, the brackets are
$$
\eqalign{
[P_\xi,P_\eta]&=P_{L_\xi\eta}\;,\cr
[P_\xi,Q_\chi]&=Q_{L_\xi\chi}\;,\cr
[Q_\chi,Q_\psi]&=P_{\star(\chi\wedge\psi)}\;.
}\eqn
$$
The Jacobi identity with three $Q$'s relies on the identity
$L_{\star(\gamma\wedge\gamma)}\gamma=0$ for a bosonic closed $2$-form $\gamma$,
which when written out in
components leads to antisymmetrisation in $6$ indices.

From a derivative expansion of the parameters, we get the level decomposition.
At level $2-2i$, there are generators $P^{a_1\ldots a_i}{}_b$,
symmetric in $(a_1\ldots a_i)$  with
$P^{a_1\ldots a_{i-1}b}{}_b=0$, \ie, in $(100i)$.
At level $1-2i$, there is $Q^{a_1\ldots a_i,bc}$, symmetric in
$(a_1\ldots a_i)$ and antisymmetric in $[bc]$, with
$Q^{a_1\ldots a_{i-1}[a,bc]}=0$, \ie, in $(001i)$.

We observe that the supermultiplet of
Section \ScalarFieldCohomologySection, related to $\BB(E_4)$ (or
$S(E_4)$) in Section \KoszulSection, is in fact the adjoint module
of $E(5,10)$.
The theorem, stating a relation between $E(5,10)$ and
$\BB(E_4)$, then immediately follows:

{\narrower\noindent\underbar{\it Theorem:} The part of $\BB(E_4)$ at
levels $\geq3$ is freely 
generated by the coadjoint module of $E(5,10)$, with the lowest
states assigned to level $3$. \smallskip}

If it holds that $S_+(E_4)\simeq\BB_+(E_4)$ (see the discussion in
Section \SEFourSection), the theorem applies to the part of
$S(E_4)$ at levels $\geq3$. This is precisely half of $S(E_4)$, and 
the theorem then gives complete information concerning the
$SL(5)$ modules in $S(E_4)$ at all levels, since the invariant
quadratic form relates the remaining half as
$R_{2-i}=\overline{R_{3+i}}$.

One may also understand $E(5,10)$ in terms of Chevalley-like
generators associated to the Dynkin diagram of Figure \BEFourFigure\
and the corresponding Cartan matrix.
The disconnectedness of the Dynkin diagram of $A_1\oplus A_2$ obtained by
deleting node $2$ presents the option not to include all three
$f_{0j}$ from the definition of $S(E_4)$\foot{The r\^ole of the
generators $f_{0j}$ in a tensor hierarchy algebra $S(\hbox{\eightfrak g})$
is as the Cartan part of the adjoint of $\hbox{\eightfrak g}^-$, the Lie algebra
obtained by deleting the node(s) connected to the fermionic one. 
The list is typically redundant. Here, by including either
$f_{03}$ or $f_{04}$, the other one can be obtained by the action of
$e_{3,4}$ and $f_{3,4}$.}.
We call such a superalgebra a restricted tensor hierarchy algebra.
By including
$f_{0j'}$, $j'=3,4$, but not $f_{01}$, one obtains only $(0011)$ at
level $-1$. Then, levels $\geq3$, as constructed in $S(E_4)$, become an
ideal, which follows from the observation that level $3$ is
annihilated by level $-1$, since
$(0011)\otimes(0001)\not\supset(1000)$. The result is $E(5,10)$.

The difference between the coadjoint module and and the superalgebra
it generates freely appears first at level $6$, where the symmetric
product of level $3$, \ie, $(0002)$, enters. Its dual module in
$S(E_4)$ is $(2000)$ at level $-1$, which accounts for the difference
between $S(E_4)$ and $E(5,10)$ at level $-1$.

Yet another superalgebra can be defined as a
restricted tensor hierarchy algebra by making the complementary choice
to the one leading to $E(5,10)$: keeping $f_{01}$ and omitting
$f_{03},f_{04}$. Then, only $(0002)$ enters at level $-1$.
Now, level $-1$ annihilates level $2$, since
$(2000)\otimes(1000)\not\supset(0010)$, so levels $\geq2$ form an
ideal.
In addition, level $-2$ is
empty (the generators in $(0002)$ have vanishing brackets among
themselves in $S(E_4)$).
The resulting superalgebra is the ``strange'' (finite-dimensional)
superalgebra $P(4)$ 
[\KacSuperalgebrasII]. This observation is due to Jakob Palmkvist
[\PalmkvistPrivate]. 



\acknowledgements The author would like to thank
Jakob Palmkvist and Ingmar Saberi for discussions, and Jakob Palmkvist
for reading the manuscript and pointing out mistakes and unclarities.


\refout

\end